\documentclass[10pt]{amsart}
 \usepackage{amssymb, amscd, amsmath, amsthm, epsfig, epsf, latexsym}

\newcommand{\compactlist}{\begin{list}{$\bullet$}{\setlength{\leftmargin}{1em}}}
\def\zz{{\bf Z}}
\def\ff{{\bf F}}
 \def\qq{{\bf Q}}
\def\cc{{\bf C}}
\def\rr{{\bf R}}
\def\co{\colon\thinspace}
 \def\cals{\mathcal{S}}
\def\calk{\mathcal{K}}
\def\co{\colon}
\newcommand{\fig}[2] { \includegraphics[scale=#1]{#2} }

\newtheorem{theorem}{Theorem}
\newtheorem{lemma}[theorem]{Lemma}

\newtheorem{definition}[theorem]{Definition}
\def\co{\colon\thinspace}

\numberwithin{equation}{section}

\begin{document}

\title{Knot 4--genus and the   rank of classes in   $  \boldsymbol{W}(\qq(\boldsymbol{t})) $ }

 \author{Charles Livingston }
 
\thanks{This work was supported in part by  NSF-DMS-0707078.}

 \address{Charles Livingston: Department of Mathematics, Indiana University, Bloomington, IN 47405 }
  
 \email{livingst@indiana.edu}


 \begin{abstract} The Witt rank $\rho(w)$ of a class $w$ in  the Witt group $ W(\ff)$   of a field with involution $\ff$  is the minimal rank of a representative of the class.  In the case of the Witt group of hermitian forms over the rational function field $\qq(t)$, we define an easily computed invariant $r(w)$ and prove that modulo torsion in the Witt group, $r$ determines $\rho$; more specifically,  $\rho(4w) = r(4w)$ for all $w \in W(\qq(t))$.  The need to determine the Witt rank arises naturally in the study of the 4--genus of knots;  we illustrate the application of our algebraic results to knot theoretic  problems, providing examples for which $r$ provides stronger bounds on the 4--genus of a knot than do classical signature bounds or Ozsv\'ath-Szab\'o and Rasmussen-Khovanov bounds.  \end{abstract}

\maketitle
 

 \section{Introduction.}
 
 For a knot $K \subset S^3$, the 4--genus of $K$, $g_4(K)$, is the minimum genus of a smoothly embedded surface in $B^4$ bounded by $K$.     Although the study of this invariant has been a focus of knot theoretic research for over 50 years, it remains an intractable invariant to compute; for instance, the determination of the 4--genus for knots with 10 or fewer crossings has just been recently completed, with even the computation for individual knots being the subject of papers, for instance~\cite{ka}. 
 
The depth of continuing interest in the 4--genus is indicated by the application of the deepest tools now available in low-dimensional topology:  Kronheimer-Mrowka's study of 4-dimensional gauge theory~\cite{km}, Ozsv\'ath-Szab\'o's  development of Heegaard-Floer theory~\cite{os1}, and Rasmussen's work on Khovanov homology~\cite{ra}, have each been used to establish Milnor's conjecture that for torus knots  $T_{p,q}$,   $g_4(T_{p,q}) =  (|p|-1)(|q|-1) /2$.

Work in the 1960s identified the central role of algebraically defined Witt groups to understanding the 4--genus. As we will review in a brief appendix, to each knot $K$ there is naturally associated a Witt class, $w_K \in W(\qq(t))$, the Witt group of hermitian forms over the rational function field, having involution induced by $t \to t^{-1}$.  (Here  $w_K $ is represented by the matrix $(1-t)V_K + (1-t^{-1})V_K^t$, where $V_K$ is an integer matrix associated to $K$, the {\it Seifert matrix}.) A fundamental result states $g_4(K) \ge \frac{1}{2}\rho(w_K)$, where $\rho$ is defined as follows:

\begin{definition} For a class $w \in W(\ff)$,   the rank of $w$,  $\rho(w)$, is  the minimum dimension of a square hermitian matrix representing $w$. 
\end{definition}

For a given class $w$, determining $\rho(w)$ can be very difficult and the most effective tools for bounding $\rho(w)$ are based on bounds on  signature functions associated to the class $w$.  A few of the early papers that applied signatures to the 4--genus are  \cite{levine, milnor, murasugi, taylor, tristram, trotter}.

The goal of this paper is to more closely examine the function $\rho(w)$.  We define an easily computed invariant $r(w)$ which provides stronger bounds on $\rho(w)$ than were previously known.  We then prove that $r(w)$ completely determines $\rho(w)$, modulo torsion in the Witt group.  More precisely, our main result states 
\vskip.15in
\noindent{\bf Theorem} {\it For all $w \in W(\qq(t))$, $\rho(w) \ge  r(w)$ and $\rho(4w) = r(4w)$.}
\vskip.15in

\noindent{\bf Outline}  Let $W(\ff)$ denote the  Witt group of  nonsingular hermitian bilinear forms over a field  $\ff$ with (possibly trivial) involution.  In the case of $\ff = \rr$ or $\ff = \cc$ (with involution given by conjugation), a relatively simple exercise shows that $\rho(w) = \sigma(w)$, where $\sigma$ is the signature.  For $\ff = \qq$ the situation is more complicated.  The diagonal form $w$ with diagonal $[1, -2]$ is not Witt trivial, and thus $\rho(w) = 2$,  but   $\sigma(w) = 0$; note however that since $\sigma(w) = 0$,  $w$ represents an element of order four in $W(\qq)$, and thus $\rho(4w) = 0$. More generally, for $w \in W(\qq)$, $\rho(4w) = \sigma(4w)$.  Details are presented in Section~\ref{secwq}.

  The arguments in   Section~\ref{secwq} are fairly basic,   but they illustrate the structure of the proof of our main theorem regarding $W(\qq(t))$.   In Section~\ref{sectiondefs} we will set    the notation to be used throughout the paper and define the function $r$.  We will also discuss explicit means for computing $r$.  
  In Section~\ref{secbounds} we will prove the first   part  of the main theorem: for $w \in W(\qq(t))$,  $\rho(w) \ge r(w)$.  Following this  we have a realization result, showing in Section~\ref{secrealize} that for any form $w$, there is a form $w'$ having an identical signature function  and for which $\rho(w') = r(w')$.  Finally, in Section~\ref{sector}  it is shown that a class with trivial signature function represents 4--torsion in $W(\qq(t))$.  The main theorem is an immediate consequence of this result.
  
  The paper concludes with Section~\ref{sectionnorm} which describes how $\rho$ leads naturally to a norm on $W(\qq(t)) \otimes \qq$ and then Section~ \ref{sectionapps} presenting examples of the computation of this norm, with specific applications to determining the 4--genus of low-crossing number knots.
  
\vskip.1in

\noindent{\it Acknowledgments} Thanks are due to Pat Gilmer and Neal Stoltzfus for early discussions regarding this work.   Thanks are also due to Jim Davis and Andrew Ranicki for their insights regarding the structure of the Witt group $W(\qq(t))$.  Special thanks go to Stefan Friedl for his thoughtful commentary on an earlier version of this paper. 
  

 \section{    $ \rho(4w) = \sigma(4w)$ for $w \in W(\qq)$.  }\label{secwq}
 
 The proof that $\rho(4w) = r(4w)$ for $w\in W(\qq(t))$ is in structure the same as the proof that $\rho(4w) = 4\sigma(w)$ for $w \in W(\qq)$.  The proof we give here is broken up into three mains steps corresponding to Sections~\ref{secbounds}, \ref{secrealize} and \ref{sector}.  A fourth concluding step is identical in both settings.
 
 \begin{theorem}\label{thmmain} If $w \in W(\qq)$ then $\rho(4w) = 4\sigma(w)$.
 
 \end{theorem}
 \begin{proof}   $ $
 \begin{enumerate}
 \item $\rho(w) \ge \sigma(w)$:    For $w \in W(\qq)$ this is immediate from the definition of  the signature.  In the case of $W(\qq(t))$ the corresponding proof will reduce to a careful algebraic calculation, based on the details of the  definition of $r(w)$  given in the next section.  The argument occupies Section~\ref{secbounds}.
 
 \item  If $w \in W(\qq)$, there is a class $w' \in W(\qq)$ with $\sigma(w) = \sigma(w')$ and $\rho(w') = \sigma(w')$:  The form $w'$ is simply the form represented by the identity matrix of dimension $\sigma(w)$.  In the case of $W(\qq(t))$, the construction of the appropriate form $w'$ for which $r(w) = r(w')  $  and $\rho(w') = r(w')$ is more delicate.
 
 \item If $\sigma(w) = 0$ then $4w = 0 $:   This depends on the structure of $W(\qq)$.  As described for instance in~\cite{mh}, there is a split short exact sequence $$0 \to W(\zz) \to W(\qq) \to \oplus_p  W(\ff_p),$$ where $\ff_p$ is the finite field with $p$ elements, $p$ a prime integer. The groups $W(\ff_p)$ are all 4--torsion, and $W(\zz) \cong \zz$, with the isomorphism given by the signature.  In the case of $W(\qq(t))$, the exact sequence is replaced with the sequence  $$0 \to W(\zz[t,t^{-1}]) \to W(\qq(t)) \to \oplus_{\alpha}  W(\qq(\alpha)),$$ where the $\alpha$ are all unit complex roots of symmetric irreducible rational polynomials.  For the analysis of these Witt groups, we turn to the references~\cite{con, lith, ranicki}.
 
 \item Conclusion:  Since the signatures are the same, we have $4w = 4w' \in W(\qq)$.   Certainly $\rho(4w') \ge 4\sigma(w')$, but by construction, $4w'$ has a representative of rank exactly $4\sigma(w')$.  The desired equality follows.  The argument in the case of $W(\qq(t))$ is identical.
 \end{enumerate}
 
 \end{proof}
 

  \section{Definition of $r(w)$.} \label{sectiondefs}
   
\begin{definition}  For a  nonsingular matrix $A$ with entries in $\qq(t)$ that is hermitian with respect to the involution induced by $t \to t^{-1}$,  the signature function   $\sigma_A'(t)$ is defined by  $\sigma_A'(t) =\text{signature} (A(e^{2\pi i t}))$.   This is   well-defined function on $[0,\frac{1}{2})$ except at the finite set of points that correspond to   poles    among the entries $A$. 
\end{definition}

Since $A$ is hermitian,  $\sigma'$ is symmetric about $\frac{1}{2}$; this justifies the restriction to the  interval $[0, \frac{1}{2})$.  As given in the next definition, taking averages and differences gives the Levine~\cite{levine} and Milnor~\cite{milnor} signature functions.  (According to Matumoto~\cite{matumoto} the jumps in the Levine signature function are determined by the Milnor signatures.  Notice that the factor of  $\frac{1}{2}$ in the definition of $J_\omega$ implies that $J_\omega(t) $ represents half the jump in the signature function at $t$.)

\begin{definition} If $w \in W(\qq(t))$, and $A$  is  a matrix representing the class $w$, then:\vskip.05in
 \begin{enumerate}
 
\item $\sigma_w(t) =\frac{1}{2}( \lim_{\tau \downarrow t  }\sigma'_A(\tau)  + \lim_{\tau \uparrow t  }\sigma_A'(\tau))$.\vskip.1in

\item $J_w(t) = \frac{1}{2}( \lim_{\tau \downarrow t  }\sigma'_A(\tau)  - \lim_{\tau \uparrow t  }\sigma_A'(\tau))$.\vskip.1in

\item For $t=0$ we define $\sigma_w(0) =  \lim_{\tau \downarrow 0  }\sigma'_A(\tau) $ and $J_w(0) =0$.
 \end{enumerate}

\end{definition}

An elementary argument shows that these functions are well defined; that is, $\sigma_w(t)$ and $J_w(t)$ depend only on the class in $W(\qq(t))$ represented by the matrix $A$.\vskip.1in

\begin{figure}[h]
\fig{1.3}{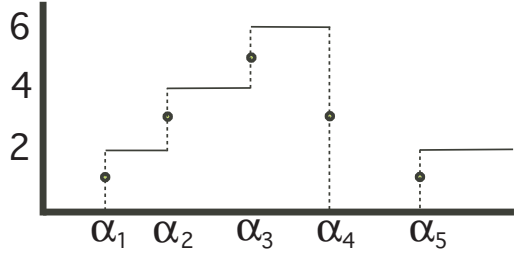}
\caption{Example of a signature function}\label{siggraph}
\end{figure}

 \noindent{\bf Example.}  Figure~\ref{siggraph} illustrates a possible signature function on the interval $[0, \frac{1}{2})$.  We will   construct a class with signature function having such a graph later, being more specific about the points $\alpha_i$.  For the specific matrix used, the values of the signatures at the discontinuities will not be known, but upon averaging, the values will be as shown in the figure.  In particular, the values of the jumps at the five discontinuities are $[1, 1, 1, -3, 1]$.

  To define the function $r \co W(\qq(t)) \to \zz_{\ge 0}$ we need to focus on the set of discontinuities of the signature function  and include in that set the value $t = 0$ for technical reasons.  That set has a natural decomposition, indexed by symmetric irreducible rational polynomials.  Throughout this paper polynomials will be Laurent polynomials $p(t) \in \qq[t, t^{-1}]$ and symmetric means $p(t) = p(t^{-1})$.  We will view these polynomials as defined on the unit complex circle.
  
\begin{definition} For  $w \in W(\qq(t))$, let $T_w$ denote  the (finite) set of discontinuities of $\sigma_w(t)$.   For each $t \in T_w$,  $e^{2 \pi i t}$ 
  is an algebraic number with rational irreducible polynomial.  For each symmetric irreducible polynomial $\delta$ 
  let $ T_{w, \delta} $ be the set $\{ t \in T_w\ | \  \delta(e^{2 \pi i t}) = 0\}$. 
   \end{definition}
 
 Given this we can define $r$.
 
 \begin{definition}  Let $w  \in W(\qq(t))$.
 
 \begin{enumerate}
 \item For each $\delta$ set $r_\delta (w) =    \max _{t \in T_{w,\delta}}   \{| \sigma_w(t)| \}+  \max_{t \in T_{w,\delta}} \{|J_w(t)| \}.$\vskip.05in
 \item $r(w) =  \max_\delta \{r_\delta(w), \sigma_\omega(0)\}.$
 \end{enumerate}
 
 \end{definition}

 \subsection{Example} Consider the graph of the signature function illustrated in Figure~\ref{siggraph}.  We assume that the five jumps occur at the roots of the same irreducible symmetric polynomial.  (We will see that as part of a general realization result, Theorem~\ref{thmrealize},   such an example does occur.  For example, the $\alpha$ can be primitive $22$--roots of units, the zeroes of  the  cyclotomic polynomial $\phi_{22}(t)$.  This is the Alexander polynomial of the torus knot $T_{2,11}$, although this torus knot does not have this signature function.)    In this example, the values of the signatures (at the $\alpha_i$ along with the value at 0) are $[0, 1,3, 5, 3, 1]$ and  the values of the jumps are $[1, 1, 1, -3, 1]$. 
 
The value of $r$ for this example will be   $$ \max\ \{1,3, 5, 3, 1\}  + \max  \{1, 1, 1,  3, 1\} ) =      5+3   = 8.$$  Notice that this is greater  than the maximum absolute value of the signature function, which is 6.
\vskip.1in

\noindent{\bf Computations}  For any hermitian matrix $A$ representing a class in $W(\qq(t))$, standard mathematical computer packages can be used to diagonalize $A$ and to arrange that the diagonal entries are Laurent polynomials.  Factoring these diagonal entries and removing factors of the form $f(t)f(t^{-1})$ ensures that these diagonal entries have factorizations as $ \delta_1(t) \cdots \delta_n(t)$, so that the $\delta_i$ are distinct symmetric polynomials   ($\delta_i(t^{-1}) =   \delta_i(t) $) with exponent one.   By symmetry, the values of the $\delta_i$ at points on the unit circle are real, and thus the signs can be determined (that is, the numerical approximation of $\delta_i(e^{i \theta})$ will be given as $a + \epsilon i$ for some small $\epsilon$, which can be ignored in the determination of the sign.)  With this, the signature function can be approximated  with necessary accuracy.  At this same time, the roots of the $\delta_i$ on the unit circle will be identified.  These computation are sufficient to completely determine the value of $r([A])$.


 \section{$r(w)$ bounds  $\rho(w)$ for $w \in W(\qq(t))$.}\label{secbounds}
 
 \begin{theorem} For any matrix representative $A$ of a class $w \in W(\qq(t))$, dim($A$) $ \ge r(w)$.
 \end{theorem}

\begin{proof}   To simplify notation, we will let $\sigma_t = \sigma_w(t)$, $\sigma_t^{+} = \lim_{\tau \downarrow t} \sigma_w(\tau)$, and  $\sigma_t^{-} = \lim_{\tau \uparrow    t} \sigma_w(\tau)$.  (For instance, for the signature function in Figure~\ref{siggraph}; $\sigma_{\alpha_2}^- = 2$ and $\sigma_{\alpha_2}^+ = 4$.)

Since $A$ is a hermitian matrix with entries in $\qq(t)$, we can diagonalize $A$,  clear denominators, and remove square factors and factors of the form $f(t)f(t^{-1})$ in the diagonal entries.  Thus, there is a diagonalization where each diagonal entry factors as the product of distinct symmetric irreducible rational polynomials.  Let  $D$   be one such diagonalization.  Note that the discontinuities of the signature function can occur only at roots of the diagonal factors.
 
  Let $T_w = \{  \alpha_1, \ldots , \alpha_k \}$ be the set of discontinuity points for the signature function on $[0,\frac{1}{2})$.

Let $\delta$ be an irreducible symmetric polynomial  such that  $\sigma_w(t)$ has a discontinuity at $t_0$ and $\delta(e^{2\pi i t_0}) = 0$.   By reordering, we can assume that $D$ has  diagonal $$[f_1 \delta, f_2 \delta, \ldots , f_m \delta,  g_1, \ldots , g_n],$$ where each $f_i$ and $g_i$ is a product of distinct irreducible symmetric polynomials, none of which are   $\delta$.

Let $\alpha  \in T_{w, \delta}$.   Evaluating the diagonal at $e^{2\pi i \alpha^-}$, where $\alpha^-$ is a number close to but smaller than $\alpha$,  we denote the count of   positive entries in $[f_1 \delta, f_2 \delta, \ldots , f_m \delta]$ to be $m_+$ and the number of negative entries to be $m_-$. Similarly,    we denote the number of positive entries in  $[g_1, \ldots , g_n]$   by   $n_+$ and the number of negative entries by  $n_-$. 

Here are some elementary calculations: 
\begin{itemize}
\item  $m_+ +   m_- = m $
\item  $n_+ + n_- = n $
\item $\sigma^-_\alpha = m_+ + n_+ - m_-  - n_-$
\end{itemize}

If we switch from $\alpha_\alpha^-$ to $\alpha_\alpha^+$ the only change in signs occurs because of the change in the sign of $\delta$ at $e^{2\pi i \alpha}$, and thus   the signs of all the diagonal entries with $\delta$ factors  change, so that we have $m_-$ positive entries and $m_+$ negative entries.  It follows, using a little arithmetic for the second calculation, that  
\begin{itemize}
\item  $\sigma^+_\alpha = m_- + n_+ - m_+  - n_-$
\item  $J_a = (m_{-} - m_+) $
\end{itemize}
Note that $J_a = m \mod 2$.
The following inequalities are verified by simply substituting for $m, J_a, n, $ and $\sigma_a^-$ in each:
\begin{itemize}

\item  $  m_\pm = \frac{1}{2} m \mp \frac{1}{2}J_a   $
\item  $n_\pm = \frac{1}{2} n \pm (\frac{1}{2} \sigma^-_a +\frac{1}{2} J_a)$
\end{itemize}
Given that $m_\pm $ and $n_\pm$ are nonnegative, we have

\begin{itemize}
\item  $ \frac{1}{2} m\ge |\frac{1}{2}J_a| $
\item  $   \frac{1}{2} n \ge  |\frac{1}{2} \sigma^-_a +\frac{1}{2} J_a|  = \frac{1}{2} |\sigma_a| $
\end{itemize}

These equations hold at each $\alpha \in T_{w, \delta}$, so, multiplying by 2 and taking the  maximums we find

\begin{itemize}
\item  $   m\ge  \max_{\alpha \in T_{w, \delta}}  | J_\alpha|$
\item  $    n \ge \max_{\alpha \in T_{w, \delta}} | \sigma_\alpha|   $
\end{itemize}

This proves the theorem, except we have not dealt yet with the signature at 0.  But clearly $\rho(w) \ge \sigma_w(0)$, the signature of $A$ evaluated near $1$ (that is, ($\sigma(0)$), so the proof is complete.
\end{proof}

\vskip.1in


\section{Realization result.}\label{secrealize}

We now want to show that every step function $s(t)$ satisfying certain criteria occurs as the signature function for some class $w \in W(\qq(t))$ and for that class, $\rho(w) = r(w)$.

\begin{definition} For a step function $s(t)$, let     $J_s(t) =  \frac{1}{2}( \lim_{\tau \downarrow t  } s(\tau)  - \lim_{\tau \uparrow t  } s(\tau))  $.\vskip.05in

\end{definition}

\begin{definition} Let $\cals$ denote the set of integer valued step functions defined on $[0,\frac{1}{2})$ satisfying the following list of conditions.  If $s \in \cals$, then:\vskip.05in

\begin{enumerate}
\item The set of discontinuities of $s$ is finite  and $s$ is continuous at $t = 0$.\vskip.05in 
\item  For all $t$, $s(t)   =\frac{1}{2}( \lim_{\tau \downarrow t  } s(\tau)  + \lim_{\tau \uparrow t  } s(\tau))$.\vskip.05in 

\item For all $t$, $ J_s(t)  \in \zz$. \vskip.05in 

\item If  $ J_s(t)  \ne 0$ then $e^{2 \pi i t}$ is the root of an irreducible symmetric rational polynomial.

\item If $\alpha_1$ and $\alpha_2$ satisfy $\delta(e^{2 \pi i \alpha_i})=0$ for some symmetric irreducible rational polynomial $\delta$, then $J_s(\alpha_1) \equiv J_s(\alpha_2) \mod 2$.

\end{enumerate}

\end{definition}

The definitions   in the previous section were given purely in terms of the signature function of a class $w \in W(\qq(t)$, so the definition extends to  $\cals$ as now described.
\begin{definition} Let $s \in S$:\vskip.05in
\begin{itemize}

\item $T_s = \{ t \ | \  J_s(t) \ne 0 \}$.

\item  For an irreducible symmetric polynomial $\delta$, $ T_{s,\delta} = \{ t\in T_s\  |\  \delta(e^{2 \pi i t}) = 0\}$.\vskip.05in

 \item For each $\delta$, $r_\delta (s) =    \max_{t \in T_{s,\delta}}   \{|  s(t)| \}+  \max_{t \in T_{s,\delta}} \{|J_s(t)| \}.$\vskip.05in
 
 \item $r(s) = \max_\delta \{r_\delta(s), s(0)\}.$
 \end{itemize}
  
\end{definition}

\begin{lemma} For all $s \in \cals$, $r(s) \equiv s(0) \mod 2$.

\end{lemma}

\begin{proof}  Since $J_s(t) \in \zz$  and $2J_s(\alpha)$ is the jump in the signature function at each discontinuity $\alpha$, $s(t) \equiv s(0) \mod 2$ if $t$ is not a point of discontinuity.  Thus, at each discontinuity $\alpha$, $s(\alpha) + J_s(\alpha) = s(\alpha)^+ \equiv s(0) \mod 2$.   For each $\delta$, for all $\alpha_1, \alpha_2 \in T_{s,\delta}$ we have $J_s(\alpha_1) \equiv J_s(\alpha_2) \mod 2$.  It follows that $s(\alpha_1) \equiv s(\alpha_2) \mod 2$.

For a fixed $\delta$, $r_\delta (s) =  \max_{t \in T_{s,\delta}}   \{| \sigma_s(t)| \} + \max_{t \in T_{s,\delta}} \{|J_w(t)| \},$ and we have now seen that mod 2, all the terms in the set of values over which the maxima are being taken are equal.  Thus, if $\alpha \in T_{s,\delta}$,  $r_\delta (s) \equiv     \sigma_s(\alpha)  +   J_s(\alpha) \mod 2$.   We have already seen that this sum equals $s(0)$, modulo 2.  

Finally, since $r(s) = \max_\delta \{r_\delta(s), s(0)\}$ and each of these elements equal $s(0)$ modulo 2, the  maximum also equals $s(0)$ modulo 2.

\end{proof}

\begin{theorem}\label{thmrealize} Suppose that $s \in \cals$.   There exists a hermitian matrix $A$ of rank $r (s)$  having signature function $s$.
  \end{theorem}
  
  \begin{proof}
  To construct $A$ we begin with the diagonal matrix $D_0$ of rank $r(s)$ in which $\delta_i$ appears as a factor of exponent one of the first  $\max_{\alpha \in T_{s, \delta_i}}  |J_s(\alpha)|$ entries.  If this condition taken over all $i$ does not specify all the entries of $D_0 $, we make the remaining entries all 1.  
 
 Next, we change the sign of some of the diagonal entries  to form $D_1$ so that the signature near 1 is $s(0)$.  This is possible, since by the previous lemma, $r(s) = s(0) \mod 2$.
 
 To continue the modification, we must introduce a family of polynomials,  $q_\theta(t) = t^{-1} - 2\cos(2 \pi \theta) +  t$,  $0 < 
\theta < \frac{1}{2}$.  
 For a dense set of $\theta$, $q_\theta(t)$ is
  a rational polynomial having its only root on the upper half circle at $e^{2 \pi \theta i}$.   For $t $ close to 0,   $q_\theta(e^{2 \pi t i})$ is positive 
  and for  $t$ close to $\frac{1}{2}$, $q_\theta(e^{2 \pi t i})$ is negative. 

If some of the diagonal entries of $D_1$ are multiplied by $q_\theta(t)$, the value of the signature is unchanged for $t < \theta$. The signature can change for values of $t > \theta$.  However, jumps can continue to appear only at the roots of the $\delta_i$, as well as, possibly, $\theta$.  The construction of the desired form consists of making such modifications to create a form with  signature function $s(t)$.

Suppose that $e^{2 \pi \alpha i }$ is a root of $\delta$, one of the $\delta_i$, corresponding to a  nontrivial jump.  Suppose also that the form $D_i$ has been constructed so that its signature function agrees with $s$ for all $t < \alpha$.  We want to alter $D_i$, building $D_{i+1}$, so that its signature function is unchanged for $t < \alpha$ and has the same jump at $\alpha$ as $s(t)$; that is,  $J_\alpha$. Pick a $\theta < \alpha$  with $\alpha - \theta$ small.  

Suppose that the first $m$ entries of $D_i$ are the ones divisible by $\delta$, and that the number of remaining entries is $n$.  If the desired form  $D_{i+1}$  is to have a jump of $2 J_\alpha$ at $\alpha$, then (when evaluated at a point $\alpha^-$ close to but less than $\alpha$) the number of positive and negative entries in $D_{i+1}$ among the first $m$ diagonal entries must be $m_+ = \frac{1}{2}m - \frac{1}{2} J_\alpha$ and $m_- =  \frac{1}{2}m + \frac{1}{2} J_\alpha$.  Similarly, if the signature to the left of $\alpha$ is to be $s_\alpha^-$, we must have the number of positive and negative entries among the last $n$ diagonal entries of $ D_{i+1}$ be   $n_+ =  \frac{1}{2}n + \frac{1}{2}s_\alpha^- + \frac{1}{2} J_\alpha$ and $n_- =  \frac{1}{2}n - \frac{1}{2}s_\alpha^- - \frac{1}{2} J_\alpha$.  (Recall that the jump is determined by the first $m$ entries, since only those change sign as $t$ increases near $\alpha$.) 

The desired sign distribution of the diagonal can be achieved by multiplying   some of the diagonal entries by $q_\theta(t)$.   The only concern is that each of the numbers  $m_+ = \frac{1}{2}m - \frac{1}{2} J_\alpha$, $m_- =  \frac{1}{2}m + \frac{1}{2} J_\alpha$, $n_+ =  \frac{1}{2}n + \frac{1}{2}s_\alpha^- + \frac{1}{2} J_\alpha$, and  $n_- =  \frac{1}{2}n - \frac{1}{2}s_\alpha^- - \frac{1}{2} J_\alpha$, must   be nonnegative.  This will be the case as long as $m \ge | \frac{1}{2}J_\alpha|$ and $n \ge |s_\alpha^- + \frac{1}{2} J_\alpha|$, which is ensured by our initial choice of the dimension of $D_0$ to be $r(s)$.  (Note that in the definition of $r(s)$  one of the two  maximum is over the numbers   $s(\alpha) $, and  $s(\alpha) = s_\alpha^- + J_\alpha$.) Observe that since $s_\alpha^-$ is unchanged, no jump has been introduced at $\theta$.
  
  \end{proof}


\section{$\rho(4w) = r(4w)$.}\label{sector}

Here we prove the main theorem:

\begin{theorem} For $w \in W(\qq(t))$, $\rho(4w) = r(4w)$.

\end{theorem}

\begin{proof} For the given form $w$, we apply Theorem~\ref{thmrealize} to the function $s = \sigma_w(t)$ to find find a hermitian matrix $A$ of dimension $r(w)$.  If we denote the class represented by $A$ in $W(\qq(t))$ by $w'$, then $\sigma_{w'}(t) = \sigma_w(t)$.  

It follows that $\sigma_{w \oplus - w'}(t) = 0$.  Lemma~\ref{lemmasig} below then shows that $w \oplus w'$ represents and element of order 1, 2, or 4, in
  $W(\qq(t))$.  Thus, $4w = 4w' \in  W(\qq(t))$.  
  Since $w'$ is constructed to have a representative of rank $r(w')$, clearly $\rho(4w') \le 4 r(w')$.  On the
 other hand, if follows immediately from the definition of $r$ that $r(nw) = n r(w)$ for any $w \in W(\qq(t))$.  Thus we have  $\rho(4w') \ge r(4w') = 4r(w')$.  The proof of the theorem is complete, given the next lemma.

\end{proof}

\begin{lemma}\label{lemmasig}  For a class $w \in W(\qq(t))$, if $\sigma_w(t) = 0$, then $4w = 0 \in W(\qq(t))$.

\end{lemma}
  
  \begin{proof}
  Background for the structure of  Witt groups is contained~\cite{mh} for symmetric bilinear forms.  The specifics in the case of hermitian forms are contained in~\cite{lith}.  A more complete description is in Ranicki's book~\cite{ranicki}, with the details of the structure of hermitian forms over number rings presented in~\cite{con}.  
  
  For each symmetric irreducible $\delta \in \qq[t,t^{-1}]$ there is a homomorphism $$\partial_\delta \co W(\qq(t)) \to W(\qq[t,t^{-1}]/\left<\delta(t)\right>),$$ defined as follows.  If $w \in W(\qq(t))$ is represented by a diagonal matrix $A$ with diagonal entries $[\delta f_1, \ldots , \delta f_m , g_1, \ldots , g_n]$ where each $f_i$ and $g_i$ is a symmetric irreducible polynomial   prime to $\delta$, then $\partial(w) = [f_1, \ldots , f_m]$.   This induces split  exact sequence $$ 0  \to  W(\qq[t,t^{-1}]) \to W(\qq(t)) \to \oplus_\delta W(\qq[t,t^{-1}]/\left<\delta(t)\right>).$$
  
  According to~\cite{ranicki} (see also \cite{lith} for an elementary argument)  the inclusion $W(\qq) \to   W(\qq[t,t^{-1}])$ is an isomorphism, so the previous sequence can be rewritten as $$ 0  \to  W(\qq) \to W(\qq(t)) \to \oplus_\delta W(\qq[t,t^{-1}]/\left<\delta(t)\right>).$$ 
  
  The field $\qq[t,t^{-1}]/\delta(t)$ is an algebraic extension of $\qq$, $\qq(\alpha)$, a field with involution given by $\alpha \to \alpha^{-1}$.  Denoting by $F$ the fixed field of the involution, an element in the Witt group of hermitian forms over this $\qq(\alpha)$ is of finite order (actually 4--torsion)  if and only for all complex embeddings of $\qq(\alpha)$ that restrict to a real embedding of $F$,  the signature of the corresponding complex hermitian form is 0.  This condition on the embedding implies that $\alpha$ maps to a unit complex root of $\delta$.
  
  Since the jump function for $w$ is 0 at all unit roots of $\delta$, it is then clear that the signature of $\partial_\delta(w)$ is also 0.  It follows that $w$ maps to an element of finite order in  $ \oplus_\delta W(\qq[t,t^{-1}]/\left<\delta(t)\right>).$ In particular, 
 $w \in W(\qq[t, t^{-1}]) \cong W(\qq)$. 
But any element in the image of $W(\qq)$ has constant
 signature function, and in our case this implies that $4w$ is represented by a class in $W(\qq)$ with 0 signature.   But $W(\qq) \cong \zz \oplus T$, where 
 $T$ is the torsion subgroup of $W(\qq)$ and satisfies $4T = 0$.  It now follows that as desired $4w = 0 \in W(\qq(t))$, using the fact the exact sequence is split exact.
    \end{proof}


 \section{$\rho$ as a norm on $W(\qq(t))$.}\label{sectionnorm}

 In order to compare $r$ as a bound on $\rho$ with bounds based on the maximum of the signature function, we want to view these functions as norms on a vector space.

 The function $r$ on $W(\qq(t))$ is multiplicative:  $r(nw) = nr(w)$ for $n \in \zz$.  Thus $r$ induces  a well-defined rational valued function on $W_\qq(\qq(t))= W(\qq(t)) \otimes \qq$.  The same is not true for $\rho$ since it can be nonzero on torsion elements in $W(\qq(t))$, but we can define a stable version of $\rho$ by $\rho_s(w) = \frac{1}{4}\rho(4w)$.  It follows then by Theorem~\ref{thmmain}  that  $\rho_s(w) = r(w)$, so $\rho_s$ also determines a well-defined rational valued function on $W_\qq(\qq(t))$.

If we define $s(w) = \max(|\sigma_w(t)|)$, then $s$ also defines a function on $W_\qq(\qq(t))$.

 \begin{theorem} For all $w \in W_\qq(\qq(t))$, $\rho_s(w) \ge s(w)$.

 \end{theorem}

   \begin{proof} Suppose the the maximum value of $|\sigma_w(t))|$ occurs at $t_0$ and $\alpha_0$ is the largest value of a discontinuity that is less than $t_0$. Then $ \sigma_w(t_0)   \le \sigma_w(\alpha_0) + J_w(\alpha_0)$.  The conclusion now follows from the definition of $ r(w)$ (which equals $\rho_s(w)$).
   \end{proof}

  Recall that a norm on a vector space $V$ is a  function $\nu$ satisfying $\nu(v ) \ge 0 $ for all $v \in V$, $\nu(v) = 0 $ if and only if $v = 0$, and $\nu(v + w) \le \nu(v) + \nu (w)$ for all $v$ and $w$.  An immediate consequence of Lemma~\ref{lemmasig} is the following.

  \begin{theorem} Both $\rho_s$ and $s$ are norms on $W_\qq(\qq(t))$.

  \end{theorem}

   \begin{definition}If $\nu$ is a norm on a vector space $V$, the unit ball of $\nu$ is defined by $B_\nu = \{ v \in V \  |\  \nu(v) \le 1\}$.

   \end{definition}
 
  
\section{Knot theoretic application; an example contrasting $\rho_s$ and $s$.}\label{sectionapps}

Here we illustrate the strength of $r$ over basic signature bounds in determining the rank of a Witt class.  We begin with a specific class, $w_1 \oplus w_2$ defined below.  We then expand on this to consider all linear combinations $xw_1 \oplus yw_2$.

\subsection{Construction and results for $w_1 \oplus w_2$} Let $\delta_6(t) = t^{-1}  -1 + t$ and let $\delta_{10} =t^{-2} - t^{-1}  +1 - t +t^2$.  These are the sixth and tenth cyclotomic polynomials, having roots at $e^{2\pi i t}$ for $t = \frac{1}{6}$ and $t = \frac{1}{10}, \frac{3}{10}$, respectively, on $[0,\frac{1}{2})$.

\begin{figure}[h]
\fig{.8}{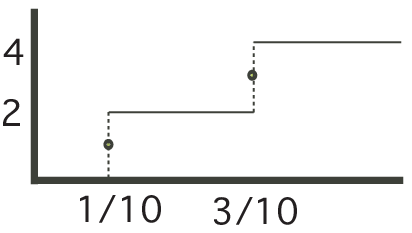}\hskip.2in\fig{.8}{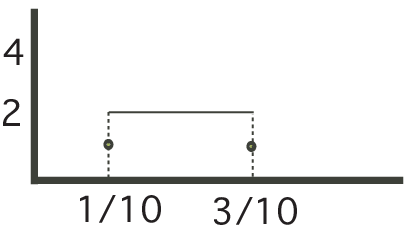}\hskip.2in\fig{.8}{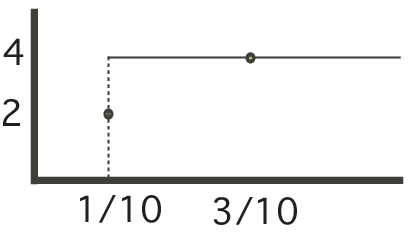}
\caption{Signature functions for $w_1$, $w_2$, and $w_1 \oplus  w_2$.}\label{fig25}
\end{figure}

We let $w_1$ be the class in $W(\qq(t))$ with diagonal representative $[- \delta_{10} \delta_6,  -\delta_6 , 1, 1] $ and let $w_2$ be the class with diagonal  $[- \delta_{10},   1] $.  The graphs of the signature functions of $w_1, w_2$, and $w_1 \oplus  w_2$ are illustrated in Figure~\ref{fig25}.  These signature functions occur for the knots $-5_1$, $10_{132}$, and $-5_1 \oplus 10_{132}$.  (The choice of signs simplifies some of the calculations that follow.) 

 The maximum  absolute values of the signature for these three forms are seen to be $s(w_1) = 4$, $s(w_2) = 2$, and $s(w_1 \oplus w_2) = 4$.  In the first two cases we have the same result for $r$:  $r(w_1) = 4$ and $r(w_2) = 2$.  However, for $w_1 \oplus w_2$ we have the set of jumps at the tenth roots of unity are $\{2, 0\}$ and the signatures at the tenth roots of unity are $\{2, 4\}$.  Thus, the sum of the two maximum is $r(w_1 \oplus w_2) = 6$.
 
 These calculations lead to the following theorem, where the knots $5_1$ and $10_{132}$  are as found in the tables at~\cite{cl}.
 
 \begin{theorem} The Witt rank of $w_1 \oplus w_2$ is $\rho(w_1 \oplus w_2) = 6$.  In particular, the knot $-5_1 \# 10_{132}$ has 4--genus 3.
 \end{theorem}
\begin{proof}  The algebraic statements are demonstrated in the discussion preceding the statement of the theorem.  For the geometric result it follows from the algebra that  $g_4(-5_1 \# 10_{132}) \ge 3$.  But it is known (for example,~see~\cite{cl}) that $g_4(5_1) = 2$ and $g_4(10_{132}) = 1$, so $g_4(-5_1 \# 10_{132}) \le 3$.
\end{proof}

\vskip.1in
\noindent{\bf Comment} This topological result can be obtained by using   Ozsv\'ath-Szab\'o invariants~\cite{os1}  or Khovanov-Rasmussen invariants~\cite{ra}, which apply only in the smooth category.  In the topological category, neither the Murasugi nor the Tristram-Levine signatures~\cite{levine, murasugi, tristram}  can give this genus bound.

\subsection{$\rho_s$ and $s$ on the span of $w_1$ and $w_2$ in $W_\qq(\qq(t))$.}  

We now compute and compare the values of $\rho_s = r$ and $s$ on the span of $w_1$ and $w_2$ in $W_\qq(\qq(t))$.  Both are determined by their unit balls.

\begin{figure}[h] 
\fig{.6}{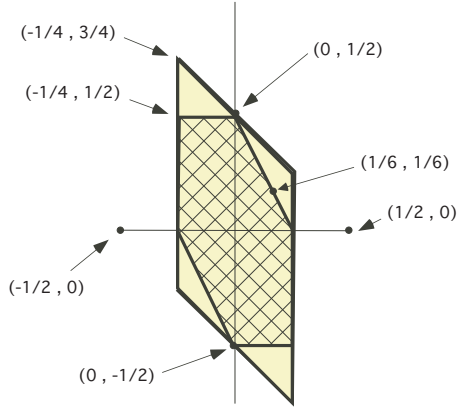}
\caption{Unit balls for $s$ and $\rho_s$.}\label{balls}
\end{figure}

The value of $s(xw_1 + yw_2)$ for $x, y \in \qq$ is given by $$s(xw_1 + yw_2)= \max\{ |2x + 2y| , |4x|\}.$$   For the value of $r$ we sum the  the maximum absolute value of the signature  at the points $t=\frac{1}{10} $ and $t=\frac{3}{10}$, and maximum for the jump function at those two points.  The result is:
$$r( xw_1 + yw_2)= \max\{ |x + y| , |3x +y|\} + \max\{ |x + y| , |x -y|\} .$$  The unit balls for these norms are drawn in Figure~\ref{balls}; the larger region represents the $s$ ball, and the smaller hatched region is the $\rho_s$ ball.

In this figure, we see that the point $(\frac{1}{4},\frac{1}{4})$ is in the unit $s$ ball; as seen earlier, the $s(w_1 \oplus w_2) = 4$.  Also, as we computed, $(\frac{1}{4},\frac{1}{4})$ is not in the unit $\rho_s$ ball, but $(\frac{1}{6},\frac{1}{6})$ is, since $\rho_s(w_1 \oplus w_2) = 6.$

\begin{figure}[h]
\fig{1.1}{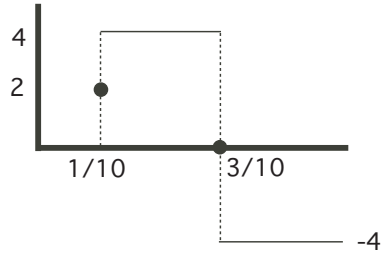}
\caption{Signature function for $-w_1 \oplus 3 w_2$.}\label{signgraph2}
\end{figure}

Another interesting point in the diagram is $(-\frac{1}{4}, \frac{3}{4})$.  The graph of the signature function of $-w_1 \oplus 3w_2$ is illustrated in Figure~\ref{signgraph2}.  The maximum absolute value of the signature is $4$, but the value of $\rho_s$ is $\max\{ 2 , 0\} + \max\{ 2, |-4|\} = 6$.

\vskip.1in

\noindent{\bf Knot theoretic comment}  It follows from this calculation that  $g_4(5_1 + 3(10_{132})) \ge 3$.   A straightforward knot theoretic exercise in fact shows shows that $g_4(5_1 \oplus 3(10_{132}) )\le 3$.  This presents another example for which $r$ detects the 4--genus of a knot, but signatures do not.  More interesting, in this case the Ozsv\'ath-Szab\'o and Rasmussen-Khovanov invariants are both insufficient to determine the 4--genus; both turn out to give a lower bound of 1.
\appendix

\section{Witt class invariants of knots}

 Here we summarize the geometric background related to $W(\qq(t))$ invariants of knots.  Details can be found in such references as~\cite{rolf}. 
 
  Every smooth oriented  knot $K\subset S^3$ bounds a smoothly embedded oriented surface $F\subset S^3$.  There is a {\it Seifert} pairing $V\co H_1(F) \times H_1(F) \to \zz$ given by $V(x, y) = lk (x, i_+(y))$, where $i_+$ is the map $F \to S^3 - F$ given by pushing off in the positive direction and $lk$ is the linking number.  A simple observation is that the intersection number of classes $x, y \in H_1(F)$ is given by $V(x,y) - V(y,x)$.  In particular, any matrix representation of $V$ has determinant $\pm 1$.
 
 Suppose the genus of $F$ is $n_1$, and $K = \partial G$, where $G$ is properly embedded in $B^4$ and is of genus $n_2$.  Then $F \cup G$ is a closed surface in $B^4$, and it bounds an embedded 3--manifold $M\subset B^4$.  An argument using Poincar\'e duality shows that the kernel $\calk$ of the inclusion $H_1(F \cup G, \qq) \to H_1(M, \qq)$ is of dimension $(n_1 + n_2)$.  Since $H_1(F)$ is a $2n_1$ dimensional subspace of $H_1(F \cup G)$, which is of dimension $2(n_1 + n_2)$, a simple linear algebra argument shows that $\calk' = \calk \cap H_1(F,\qq)$ is of dimension at least $n_1 - n_2$.
 
 Another simple geometric argument implies that $V$ vanishes on the subspace of $\calk' \subset H_1(F, \qq)$.  If we now write $V$ for a matrix representation of the Seifert pairing, the form $(1-t)V + (1 -t^{-1})V $ defines a Hermitian pairing on the rational function field.  As above, if $K$ bounds a surface of genus $n_2$ in $B^4$, then this forms vanishes on a subspace in $H_1(F,\qq)$ of dimension $(n_1 - n_2)$. Thus the form splits as a direct sum of forms, one of which is metabolic and of dimension $2(n_1 -n_2)$; the other summand is  of dimension $2 n_1   - (2(n_1 - n_2) = 2n_2$.

In summary, we see that if a knot $K$ bounds a surface of genus $g$ in $B^4$, then the Witt class of  its hermitianized Seifert form has a representative of dimension $2g$.
  


\newcommand{\etalchar}[1]{$^{#1}$}

\end{document}